\documentclass{elsart}

\usepackage{latexsym}
\usepackage{amsfonts}
\usepackage{graphicx,amssymb,cite}
\usepackage{amsmath}
\usepackage{amssymb}
\usepackage[mathscr]{eucal}
\usepackage{enumerate}

\begin{document}

\newtheorem{tm}{Theorem}[section]
\newtheorem{pp}{Proposition}[section]
\newtheorem{lm}{Lemma}[section]
\newtheorem{df}{Definition}[section]
\newtheorem{tl}{Corollary}[section]
\newtheorem{re}{Remark}[section]
\newtheorem{eap}{Example}[section]

\newcommand{\pof}{\noindent {\bf Proof} }
\newcommand{\ep}{$\quad \Box$}

\newcommand{\al}{\alpha}
\newcommand{\be}{\beta}
\newcommand{\var}{\varepsilon}
\newcommand{\la}{\lambda}
\newcommand{\de}{\delta}
\newcommand{\st}{\stackrel}

\allowdisplaybreaks

\begin{frontmatter}

\title{ A countable set derived by fuzzy set}
\author{Huan Huang}
\author{}\ead{hhuangjy@126.com (H. Huang)}
\address{ Department of Mathematics, Jimei
University, Xiamen 361021, China. }
\date{}
\maketitle

\begin{abstract}
In this paper, it shows that
for each fuzzy set $u$ on $\mathbb{R}^m$,
  the set
 $D(u)$ is at most countable.
Based on this, it modifies 
the proof of 
assertion (I) in step 2 of the sufficiency part of Theorem 4.1 in paper: Characterizations of compact sets in fuzzy sets spaces with $L_p$ metric, http://arxiv.org/abs/1509.00447.

\end{abstract}


\end{frontmatter}

\section{Introduction}

In this paper, we show that
for each fuzzy set $u$ on $\mathbb{R}^m$,
  the set
 $D(u)=\{\al\in (0,1):  [u]_\al \nsubseteq  \overline{ \{u> \al\}}  \;\}$ is at most countable.

\section{Preliminaries}

We introduced some basic concepts about fuzzy sets. For details, we refer to \cite{wu,da}.

We use $F(\mathbb{R}^m)$ to represent all
fuzzy subsets on $\mathbb{R}^m$, i.e. functions from $\mathbb{R}^m$
to $[0,1]$.
For
$u\in F(\mathbb{R}^m)$, let $[u]_{\al}$ denote the $\al$-cut of
$u$, i.e.
\[
[u]_{\al}=\begin{cases}
\{x\in \mathbb{R}^m : u(x)\geq \al \}, & \ \al\in(0,1],
\\
{\rm supp}\, u=\overline{\{x \in \mathbb{R}^m: u(x)>0\}}, & \ \al=0.
\end{cases}
\]
Let $\{u>\al\}$ denote the strong $\al$-cut of $u$, i.e.,
$\{u>\al\}= \{   x \in \mathbb{R}^m: u(x)>\al   \}$.

\section{Main results}
Let $u\in F(\mathbb{R}^m)$, $t\in \mathbb{R}^m$ and $r$ be a positive number in $\mathbb{R}$.
Define a function
$S_{u,t,r}(\cdot, \cdot): \mathbb{S}^{m-1} \times [0,1] \to \{ -\infty \} \cup \mathbb{R}$
by
$$
S_{u,t,r}(e, \al)= \left\{
                        \begin{array}{ll}
      -\infty,    &     \hbox{if} \ [u]_\al \cap \overline{ B(t, r)}  =\emptyset,
\\
\sup \{    \langle e, x-t \rangle : x\in [u]_\al \cap \overline{ B(t, r)}   \} ,  &  \hbox{if} \ [u]_\al \cap \overline{ B(t, r)}  \not=\emptyset,
                        \end{array}
                      \right.
  $$
where  $\overline{ B(t, r)}$ denote
the closed ball
$\{x\in \mathbb{R}^m: d(t,x) \leq r \}$.

We say $\alpha\in (0,1)$ is a discontinuous point of $S_{u,t,r}(e, \cdot)$
if
\\
(\romannumeral1)  $S_{u,t,r}(e, \alpha) \in \mathbb{R}$, and
\\
  (\romannumeral2) $S_{u,t,r}(e, \beta) = -\infty$ for all $\beta>\alpha$
or
$-\infty < \lim_{\beta \to \alpha+}S_{u,t,r}(e, \beta)  <   \lim_{\beta \to \alpha-}S_{u,t,r}(e, \beta) $.
\\
Denote the set of all discontinuous points of $S_{u,t,r}(e, \cdot)$ by $D_{u, t, r, e}$.
Then
$D_{u,t,r,e}$ is at most countable
because $S_{u,t,r}(e, \cdot)$ is a monotone function on $[0,1]$.

\begin{tm} \label{dutrcoun}
  Let $u$ be a fuzzy set on $\mathbb{R}^m$, $t$ be a point in $\mathbb{R}^m$, and $r$ be a positive real number.
  Then
  $D_{u,t,r}:=\bigcup_{e\in \mathbb{S}^{m-1}} D_{u,t,r,e}$
  is at most countable.
\end{tm}

\pof \ Let $\varpi$ be a countable dense subset of $\mathbb{S}^{m-1}$.
Then
\begin{equation}\label{couren}
 D_{u,t,r} =\bigcup_{e\in \varpi} D_{u,t,r,e}
\end{equation}
In fact,
suppose that $\al\in D_{u,t,r}$,
then there exists
$e\in \mathbb{S}^{m-1}$ such that
 $\alpha \in D_{u,t,r, e}  $.
Hence
$S(u,t,r)(e, \alpha)>-\infty$.
This is equivalent to
$[u]_\alpha \cap \overline{B(t,r)} \not=\emptyset.$
Therefore
\begin{equation}\label{sar}
S(u,t,r)(f, \alpha)>-\infty \   \hbox{for all} \  f\in \mathbb{S}^{m-1}.
\end{equation}
To show $\al \in \bigcup_{e\in \varpi} D_{u,t,r,e}$,
we divide the proof into two cases.

Case 1. \ $S(u,t,r)(e, \beta)=-\infty$
for all
$\beta > \alpha$.

In this case,
$[u]_\beta \cap \overline{ B(t, r)}  =\emptyset$ for any
$\beta > \alpha$,
and so
$S(u,t,r)(f, \beta)=-\infty$ when $f\in \mathbb{S}^{m-1}$ and $\beta > \alpha$.
Combined with \eqref{sar}, we know
$\al \in D_{u,t,r,f}$ for each $f\in \mathbb{S}^{m-1}$.
Thus
$\al \in \bigcup_{e\in \varpi} D_{u,t,r,e}$.

Case 2. \ $-\infty < \lim_{\beta \to \alpha+}S_{u,t,r}(e, \beta)  <   \lim_{\beta \to \alpha-}S_{u,t,r}(e, \beta) $.

In this case, there is an $\al_0 > \al $ such that $[u]_{\lambda} \cap  \overline{B(t,r)} \not= \emptyset $
 when $\lambda \in [0, \al_0]$.
Set
\begin{equation}\label{den}
  \xi : =  \lim_{\beta \to \alpha-}S_{u,t,r}(e, \beta)  -   \lim_{\beta \to \alpha+}S_{u,t,r}(e, \beta)  >0.
\end{equation}
 Notice that, for all $\beta\in [0,1]$ with $[u]_\beta \cap \overline{B(t,r)} \not= \emptyset$,
 \begin{align*}
    | S_{u,t,r}(e, \beta) &-  S_{u,t,r}(f, \beta)|    \\
&=
  | \sup \{    \langle e, x-t \rangle : x\in [u]_\beta  \cap \overline{B(t,r)}  \} - \sup \{    \langle f, x-t \rangle : x\in [u]_\beta  \cap \overline{B(t,r)}  \} |
   \\
    & \leq  \sup \{  \mid \langle  e-f,   x-t \rangle \mid : x\in [u]_\beta \cap \overline{B(t,r)} \}
    \\
    & \leq  \|  e-f \|  \cdot r,
 \end{align*}
hence, for any $\lambda \in [0, \al_0]$,
$$ \mid  S_{u,t,r}(e, \lambda) -  S_{u,t,r}(f, \lambda)\mid \leq  \|  e-f \|  \cdot r,$$
and so, combined with \eqref{den},
we know
there exists $\delta>0$ such that, for all $f\in  \mathbb{S}^{m-1} \cap B(e,\delta)$,
 $$ \lim_{\beta \to \alpha-}S_{u,t,r}(f, \beta)  -   \lim_{\beta \to \alpha+}S_{u,t,r}(f, \beta) > \xi/2, $$
 this means that $\alpha \in D_{u,t,r,f}$ when $f\in  \mathbb{S}^{m-1} \cap B(e,\delta)$.
Thus there exists $g\in \varpi$
such that
$\al\in D_{u,t,r,g}$,
i.e.,
$\al \in \bigcup_{e\in \varpi} D_{u,t,r,e}$.

Now we obtain \eqref{couren}. Since $\varpi$ is countable and $D_{u,t,r,e}$ is at most countable,
we know
$D_{u,t,r}$ is at most countable.
\ep

\begin{re}{\rm
  In the proof of Theorem \ref{dutrcoun}, in order to show $ D_{u,t,r}$ is at most countable,
it proves
that
 $D_{u,t,r}=\bigcup_{e\in \varpi} D_{u,t,r,e}$.
This kind of trick was used in the proof of
Lemma 4 in \cite{trutschnig} to show a set is at most countable.
}
\end{re}

\begin{tm} \label{fusdctn}
  Let $u$ be a fuzzy set on $\mathbb{R}^m$,
  then $D(u)=\{\al\in (0,1):  [u]_\al \nsubseteq  \overline{ \{u> \al\}}  \;\}$ is at most countable.
\end{tm}

\pof \ If $\alpha \in D(u)$, then
 there is a
$y\in \mathbb{R}^m$ such that $y\in [u]_\al$
but
$y \notin \overline{ \{u>\al\}}$.
Thus,
\begin{equation}\label{dgvn}
  d(y, \overline{ \{u>\al\}}) > \varepsilon>0.
\end{equation}
Choose a $q \in \mathbb{Q}^m = \{(z_1, z_2,\ldots, z_m) \in \mathbb{R}^m :   z_i \in \mathbb{Q}, \ i=1,2,\ldots, m \}   $ which
satisfies that
 $\|y-q\|>0$. We assure  that $\al \in D_{u,q, r}$ for some $r\in Q$ with $r\geq \|y-q\|$.

In fact, let $e= (y-q)/\|y-q\|  $. Then
\begin{equation}\label{suqrae}
  S_{u,q,r}(e, \al)\geq \langle e, y-q \rangle = \|y-q\|.
\end{equation}
If $[u]_\beta \cap \overline{B(q,r)} = \emptyset$ for   any $\beta>\alpha$,
then
$S_{u,q,r}(e, \beta)=-\infty$ for all $\beta > \alpha$,
and thus
$\al \in D_{u,q, r,e}$.

If    there exists $\beta>\alpha$ such
that
$[u]_\beta \cap \overline{B(q,r)} \not= \emptyset$.
Pick an arbitrary
$x\in \overline{\{u>\al\}} \cap \overline{B(q,r)} $,
then, by \eqref{dgvn},
\begin{gather*}
  \|x-q\| \leq r,
 \\
\|x-y\| >\varepsilon.
\end{gather*}
If $x=q$, then $ \langle  e, x-q \rangle =0$.
Suppose that $x \not= q$.
Notice that
$$ \frac{\langle y-q , x-q\rangle}{ \|y-q\| \cdot \|x-q \|}=\cos \alpha
=
 \frac{\|x-q\|^2 +  \|y-q\|^2 -\|x-y\|^2 }{2\|y-q\| \cdot   \|x-q\| },$$
where
$\al$ is the angle between two vectors
$x-q$ and $y-q$.
Thus
\begin{align*}
  \langle  e, x-q \rangle & =\frac{1}{2}
\left(     \frac{\|x-q\|^2}{ \|y-q\|} +   \|y-q\|  -  \frac{ \|x-y\|^2}{ \|y-q\|}   \right )
\\
   & \leq   \frac{1}{2} \left(    \frac{r^2}{ \|y-q\|}   +   \|y-q\|  -  \frac{ \varepsilon^2}{ \|y-q\|}   \right ),
\end{align*}
and
so
 there exists a $\delta>0$ such that for all $r \in [ \|y-q\|,  \|y-q\|+\delta  )$,
\begin{equation}\label{suqrbe}
   \langle  e, x-q \rangle   \leq    \|y-q\| -   \frac{1}{4} \frac{ \varepsilon^2}{ \|y-q\|}.
\end{equation}
Combined with \eqref{suqrae} and \eqref{suqrbe},
it then follows from the arbitrariness
of
$x\in \overline{\{u>\al\}} \cap \overline{B(q,r)} $
that
$$\lim_{\beta \to \al+} S_{u,q,r}(e,\beta) <  S_{u,q,r}(e,\alpha)$$
when
$r \in [ \|y-q\|,  \|y-q\|+\delta  )$.
This implies
that
there exists $r\in \mathbb{Q}$
such that
$\al \in D_{u,q,r,e} \subset  D_{u,q,r}$.

Now we know
$$D(u) \subseteq   \bigcup_{q\in \mathbb{Q}^m, r\in \mathbb{Q}} D_{u,q,r}.$$
It then follows from Theorem \ref{dutrcoun}
that
 $D(u)$
is at most countable. \ep

\section{An application}

In this section, we give an application of the main result.

The following is
assertion (I) in step 2 of the sufficiency part of the proof of Theorem 4.1 in paper: Characterizations of compact sets in fuzzy sets spaces with $L_p$ metric, see http://arxiv.org/abs/1509.00447.

(I)\ $P(v)$ is at most countable.

We can show this assertion in the following way.

In fact, given $\al\in P(v)$, it holds that
 $\overline{\{v>\alpha\} }\varsubsetneqq   [v]_\al$,
and hence
$\al \in D(v):=\{\gamma\in (0,1):    [v]_\gamma \nsubseteq  \overline{ \{v> \gamma\}}  \; \}$.
Thus $P(v) \subseteq D(v)  $.
By Theorem \ref{fusdctn},
$D(v) $ is at most countable.
So
$P(v)$ is at most countable.

\end{document}